\documentclass[reqno]{amsart}
\usepackage{calrsfs}
\DeclareMathAlphabet{\pazocal}{OMS}{zplm}{m}{n}
\usepackage{accents,amssymb,amscd,bbm,enumitem,graphicx,mathtools,verbatim}
\usepackage[mathscr]{euscript}
\usepackage[colorlinks]{hyperref}
\newtheorem{theorem}{Theorem}[section]
\newtheorem{corollary}[theorem]{Corollary}
\newtheorem{lemma}[theorem]{Lemma}

\theoremstyle{definition}

\newtheorem{hyp}[theorem]{Hypothesis}

\numberwithin{equation}{section}
\DeclarePairedDelimiter{\norm}{\lVert}{\rVert}

\renewcommand{\Im}{\operatorname{Im}}
\newcommand{\supp}{\operatorname{supp}}

\newcommand{\iu}{\mathrm{i}}
\newcommand{\m}{\mathrm{m}}
\newcommand{\dom}{\operatorname{dom}}
\newcommand{\ran}{\operatorname{ran}}

\newcommand{\bv}{\operatorname{BV}}
\newcommand{\bvl}{\operatorname{BV}_{\rm loc}}
\newcommand{\bvup}{\operatorname{BV}_{\uparrow}}
\newcommand{\bvdown}{\operatorname{BV}_{\downarrow}}

\newcommand{\sgn}{\operatorname{sgn}}

\newcommand{\sm}[1]{\big(\begin{smallmatrix}#1\end{smallmatrix}\big)}
\newcommand{\bsm}[1]{\begin{pmatrix}#1\end{pmatrix}}
\newcommand{\bb}[1]{{\mathbb{#1}}}
\newcommand{\br}[1]{{\breve{#1}}}
\newcommand{\ms}[1]{{\mathscr{#1}}}
\newcommand{\mc}[1]{{\mathcal{#1}}}
\newcommand{\pc}[1]{{\pazocal{#1}}}
\newcommand{\id}{\mathbbm 1}
\newcommand{\ov}{\overline}
\newcommand{\ol}{{\overline{\lambda}}}
\newcommand{\la}{\lambda}

\newcommand{\loc}{{\rm loc}}
\newcommand{\mx}{{\rm max}}
\newcommand{\mn}{{\rm min}}

\newcommand{\sas}{\;\;\text{and}\;\;}
\newcommand{\<}{\langle}
\renewcommand{\>}{\rangle}

\newcommand{\cev}[1]{\overset{{}_{\leftharpoonup}}{#1}}
\newcommand{\brx}{\xi}
\newcommand{\syc}{C}
\begin{document}
\title[Liouville transforms]{Liouville transforms for ordinary differential equations with distributional coefficients}
\author{Rudi Weikard}
\address{Department of Mathematics, University of Alabama at Birmingham, Birmingham, AL 35226-1170, USA}
\email{weikard@uab.edu}
\date{\today}



\begin{abstract}
We study Liouville transforms for the differential equation $Ju'+qu=wf$ on the real interval $(a,b)$ where $J$ is a constant, invertible, skew-hermitian matrix and $q$ and $w$ are matrices whose entries are distributions of order $0$ with $q$ Hermitian and $w$ non-negative.
In particular, we examine their use in the treatment of quasi-regular endpoints and their impact on the Weyl-Titchmarsh $\bb M$-function associated with the differential equation.
\end{abstract}
\maketitle

\section{Introduction}
In this paper we investigate Liouville transforms for the differential equation
\begin{equation}\label{de}
Ju'+qu=wf
\end{equation}
on the interval $(a,b)\subset\bb R$ where the coefficients $J$, $q$, and $w$ satisfy the following hypothesis.
\begin{hyp}\label{hyp:m}
$J$ is a constant, invertible, skew-hermitian $n\times n$-matrix and $q$ and $w$ are $n\times n$-matrices whose entries are distributions of order $0$ on $(a,b)$; $q$ is Hermitian and $w$ is nonnegative.
\end{hyp}

We also require that $f:(a,b)\to\bb C^n$ is such that the components of $wf$ are still distributions of order $0$. 
The solutions of equation \eqref{de} are then to be sought among the balanced%
\footnote{This term will be defined near the end of this introduction.}
functions $u:(a,b)\to\bb C^n$ of locally bounded variation. 
Since the (distributional) derivatives of functions of locally bounded variation are distributions of order $0$, equation \eqref{de} is a well-posed problem.

There is a one-to-one correspondence between distributions of order $0$ and local measures on $(a,b)$. 
A local measure (for our purposes) is a complex-valued, countably additive function defined on the ring $\mc B'$ of all Borel subsets of $(a,b)$ that are included in some compact subset of $(a,b)$.
Here, to be precise, countable additivity of the local measure $\mu$ means that
$\mu(A)=\sum_{n=1}^\infty \mu(A_n)$ provided that $n\mapsto A_n$ is a sequence of pairwise disjoint sets in $\mc B'$ and that $A=\bigcup_{n=1}^\infty A_n$ is also in $\mc B'$. We emphasize that, if $K$ is a compact subset of $(a,b)$, then $\{E\cap K:E\in\mc B'\}$ is a $\sigma$-algebra and the restriction of $\mu$ to this set is a complex measure.
The most important facts about local measures are collected in Appendix \ref{A:lm} but a more complete presentation may be found in 
\cite{MR5052520}.
The aforementioned correspondence between distributions of order $0$ and local measures on $(a,b)$ is now seen as follows:
If $\phi\in C^\infty_c((a,b))$ is a test function, then the map $\phi\mapsto \int_{(a,b)}\phi\mu$ is a distribution of order $0$ whenever $\mu$ is a local measure. 
That all distributions of order $0$ are obtained in this way follows from Riesz's representation theorem, see, for instance, Section 6.2 of \cite{MR5052520}. 
Another, not less important fact, is that a cumulative distribution function $M$ of a local measure $\mu$ is a function of locally bounded variation and that the Lebesgue-Stieltjes measure generated by $M$ is $\mu$, i.e., we have $\mu=dM$. 
The function $M$ gives also rise to a distribution of order $0$ via the assignment $\phi \mapsto \int_{(a,b)}\phi M\,dx$%
\footnote{Here we use the symbol $dx$ to signify Lebesgue measure.}.
We often identify, as is common, the function $M$ with the distribution $\phi \mapsto \int_{(a,b)}\phi M\,dx$. 
Its distributional derivative $M'$ is given by $\phi \mapsto \int_{(a,b)}\phi\,dM$.
We therefore use the notations $\mu$, $dM$, and $M'$ interchangeably and, accordingly do not distinguish between a local measure and its corresponding distribution of order $0$.
Given this understanding of a distribution of order $0$ we define now the product of such a distribution $r$ with a suitable function $f$ by setting
$$(rf)(\phi)=\int f\phi dR$$
when $R$ is an antiderivative of $r$; in order for the integral here to make sense it is required that $f$ is locally integrable with respect to the total variation of $dR$.

If the entries of $q$ and $w$ are measures that are locally absolutely continuous with respect to Lebesgue measure, it is common to interpret $q$ and $w$ as the corresponding Radon-Nikodym derivatives, i.e., locally integrable (with respect to Lebesgue measure) functions. 
In this case equation \eqref{de} is, of course, well studied and encompasses also differential equations of higher order when written as a system. 
One may now raise the questions whether $q$ and $w$ could be even rougher than distributions of order $0$ are. 
If they were, the solutions $u$ would also be rougher and it would become impossible to define the products $qu$ and $wu$. 

In our case, where $q$ and $w$ are matrices of distributions of order $0$, the corresponding local measure are matrix-valued. 
In particular, if $E$ is a Borel set that is included in a compact subset of $(a,b)$, then $q(E)$ is a Hermitian matrix and $w(E)$ is a nonnegative matrix.
 
The first (to the best of my knowledge) to consider a Sturm-Liouville equation with measure coefficients, a special case of the present set-up, was Krein \cite{MR0054078} in 1952 when he modeled a vibrating string.
Further contributions we are aware of were made by Kac \cite{MR0080835}, Feller \cite{MR0068082}, Mingarelli \cite{MR706255}, Gesztesy and Holden \cite{MR914699}, Kurasov \cite{MR1397901}, Kurasov and Boman\cite{MR1443392}, Savchuk and Shkalikov \cite{MR1756602}, Eckhardt et al. \cite{MR3046408}, and Eckhardt and Teschl \cite{MR3095152}.
We refer the reader to Ghatasheh and Weikard \cite{MR4047968} and Eckhardt et al.\ \cite{MR3046408} for more details on the history of the subject.
For general information about Sturm-Liouville equations consider the classic book \cite{MR2170950} by Zettl or the quite recent book \cite{MR4809879} by Gesztesy, Nichols, and Zinchenko.

The material in this paper is organized as follows.
In Section \ref{S:prelim} we present the basic facts about systems of ordinary differential equations with distributional coefficients and their associated linear relations. 
Then we define and investigate Liouville transforms in Section \ref{S:Liouville}. 
In Section \ref{S:bc} we show how Liouville transforms may be used to pose initial conditions at a quasi-regular%
\footnote{The term quasi-regular will be defined in Section \ref{S:prelim}.} endpoint.
In Section \ref{S:WTM} we recall the construction of the Weyl-Titchmarsh function $\bb M$ for a self-adjoint linear relation associated with our differential equation \eqref{de}. 
In the present generality this has been done in \cite{MR4711083} to which we refer for proofs. 
Then we investigate, in Section \ref{S:WTMLT}, the relationship of two Weyl-Titchmarsh functions $\bb M$ and $\br{\bb M}$, respectively associated with two self-adjoint linear relations $T$ and $\br T$ that are connected via a Liouville transform.
Finally, in Section \ref{S:n=2}, we treat the special case where $n=2$ and the left endpoint of the interval $(a,b)$ is quasi-regular.
It turns out, not unexpectedly, that the Weyl-Titchmarsh matrix $\bb M$ is then determined by a scalar function $m$. 

We end this introduction with a few words about notation.
The identity matrix is denoted by $\id$ while $\chi_E$ stands for the characteristic function of the set $E$. 
The space of functions of locally bounded variation defined on $(a,b)$ with values in $\bb C^{m\times n}$ is denoted by $\bvl((a,b))^{m\times n}$. 
Functions in $\bvl((a,b))^{m\times n}$ that are of bounded variation on $(a,c)$ for some (and hence all) $c\in(a,b)$ are collected in the space
$\bvdown((a,b))^{m\times n}$ while $\bvup((a,b))^{m\times n}$ stands for the collection of functions in $\bvl((a,b))^{m\times n}$ that are of bounded variation on $(c,b)$ for some $c\in(a,b)$. 
Finally, $\bv((a,b))^{m\times n}$ denotes the space of functions of bounded variation on $(a,b)$. 

Any function in $u\in\bvl((a,b))^{m\times n}$ has left- and right-hand limits at any $x\in(a,b)$. 
Thus we may introduce the left-continuous function $x\mapsto u^-(x)=\lim_{t\uparrow x} u(x)$ and the right-continuous function $x\mapsto u^+(x)=\lim_{t\downarrow x} u(x)$. 
Moreover, we let 
$$u^\#=\frac12(u^-+u^+).$$ 
A function $u$ of locally bounded variation is called \emph{balanced} if $u=u^\#$. 
The balanced functions of locally bounded variation on $(a,b)$ form a linear subspace of $\bvl((a,b))^{m\times n}$ denoted by $\bvl^\#((a,b))^{m\times n}$.
We define, of course, the spaces $\bvdown^\#((a,b))^{m\times n}$, $\bvup^\#((a,b))^{m\times n}$, and $\bv^\#((a,b))^{m\times n}$ analogously. 

An $n\times n$-matrix $S$ is called \emph{symplectic} (with respect to $J$) if $S^*JS=J$.
Since $|\det S|^2=1$ it follows that symplectic matrices are invertible.
Indeed, their inverses are again symplectic.
Moreover, products of symplectic matrices are symplectic. 
In the special case where $n=2$ and $J=\sm{0&-1\\ 1&0}$ it is the case that a matrix is symplectic if and only if it is the product of a complex number of modulus $1$ and a real matrix of determinant $1$. 
We also record the following fact about symplectic matrices.
\begin{lemma}\label{L:2symps}
If $A$ and $B$ are two matrices, symplectic with respect to $J$, then 
$$(A^*-B^*)(A^*+B^*)^{-1}=J(A+B)^{-1}(B-A)J^{-1}.$$
\end{lemma}

\begin{proof}
Let $Y=A^{-1}B$.
Then $(A^*-B^*)(A^*+B^*)^{-1}=(\id-Y^*)(\id+Y^*)^{-1}$.
Since $(\id+Y^*)^{-1}$ and $\id-Y^*$ commute we get $(A^*-B^*)(A^*+B^*)^{-1}=(\id+Y^*)^{-1}(\id-Y^*)$.
Now notice that $Y$ is also symplectic so that $Y^*=JY^{-1}J^{-1}=JB^{-1}AJ^{-1}$. 
Using this property one gets the result after a few further manipulations. 
\end{proof}

\section{Preliminaries}\label{S:prelim}
In this section we gather some basic material on systems of ordinary differential equations with locally integrable coefficients and the associated linear relations for easy reference. 

\subsection{Background on differential equations}
Since $w$ is a nonnegative measure one may define the space $\pc L^2(w)$ consisting of those functions $f:(a,b)\to\bb C^n$ for which the measure $f^*wf$ is defined and finite on $(a,b)$. 
In $\pc L^2(w)$ one defines the semi-scalar product $\<f,g\>=\int f^*wg$ and the semi-norm $\norm f=\<f,f\>^{1/2}$.
Identifying functions in $\pc L^2(w)$ whose difference have semi-norm $0$ we obtain the Hilbert space $L^2(w)$.

As mentioned above we require the solution of equation \eqref{de} to be a balanced function of locally bounded variation.  
The choice of considering only balanced functions is motivated by the following observation. 
Let $u,v\in\bvl((a,b))^{n}$ satisfy $u=\kappa u^++(1-\kappa)u^-$ and $v=\kappa v^++(1-\kappa)v^-$ for some fixed parameter $\kappa$.
Then the integration by parts formula implies that
\begin{equation}\label{ibypt}
\int_{[c,d]} (u^*dv+(du)^*v)=(u^*v)^+(d)-(u^*v)^-(c)+(2\kappa-1)\int_{[c,d]}(u^+-u^-)^*dv
\end{equation}
whenever $[c,d]\subset(a,b)$.
Thus, when $\kappa=1/2$, i.e., when $u$ and $v$ are balanced, the usual form of this formula -- essential for spectral theory -- holds. 

Given a distribution $r$ of order $0$ with antiderivative $R$ we set, for any $x\in(a,b)$,
$$\Delta_r(x)=R^+(x)-R^-(x)=dR(\{x\}).$$
Note that the set of those $x$ where $\Delta_r(x)\neq0$ is at most countable. 
With this notation we now define the matrices 
$$B_\pm(x,\la)=J\pm\frac12(\Delta_q(x)-\la\Delta_w(x)).$$
We note, for future reference, that
\begin{equation}\label{bpmstar}
B_\pm(x,\la)^*=-B_\mp(x,\ol).
\end{equation}
The differential equation $Ju'+(q-\la w)u=wf$ is satisfied if $Ju'+(q-\la w)u-wf$ is the zero-measure.
In particular, evaluating this measure at the singleton $\{x\}$ yields the necessary condition that 
$$B_+(x,\la)u^+(x)-B_-(x,\la)u^-(x)=\Delta_w(x)f(x)$$
must hold for every $x\in(a,b)$. 
This becomes $J(u^+(x)-u^-(x))=0$, if $\Delta_q(x)=\Delta_w(x)=0$, i.e., the balanced function $u$ must be continuous at such points. 

If $\Delta_q(x)$ or $\Delta_w(x)$ are different from $0$ one or both of the matrices $B_\pm(x,\la)$ may fail to be invertible. 
Then it is impossible to uniquely determine $u^+(x)$ or $u^-(x)$ even if the other one is given. 
We collect the values of $\la$ causing such trouble in the set 
$$\Lambda_x=\{\la\in\bb C: \text{$\det B_+(x,\la)=0$ or $\det B_-(x,\la)=0$}\}.$$ 
We also define 
$$\Lambda=\bigcup_{x\in(a,b)}\Lambda_x$$
and 
$$\Xi_\la=\{x\in(a,b): \text{$\det B_+(x,\la)=0$ or $\det B_-(x,\la)=0$}\}.$$ 

We emphasize that $\Lambda_x$ is empty, finite, or all of $\bb C$ and that $\Lambda_x\neq\emptyset$ can only happen for countable many values of $x$.
It follows that $\Lambda$ is either equal to $\bb C$ or else a countable set.
Moreover, on account of equation \eqref{bpmstar}, we have $\Xi_\la=\Xi_\ol$ and that each $\Lambda_x$ is symmetric with respect to the real axis. 

We introduce the following classification of the endpoints $a$ and $b$ of the interval $(a,b)$ on which our problem is posed. 
Let $Q$ and $W$ be antiderivatives of $q$ and $w$, respectively. 
These are elements of $\bvl((a,b))^{n\times n}$.
The endpoint $a$ is called \emph{regular}, if both $Q$ and $W$ are in $\bvdown((a,b))^{n\times n}$.
Equivalently, $a$ is regular if $q$ and $w$ are finite measures on $(a,c)$ for some (and hence all) $c\in(a,b)$. 
Similarly, $b$ is called a regular endpoint if $Q$ and $W$ are in $\bvup((a,b))^{n\times n}$ or if $q$ and $w$ are finite measures on $(c,b)$ for some $c\in(a,b)$.
An endpoint that is not regular is called \emph{singular}. 
Of course, this classification depends on the coefficients $q$ and $w$.

The endpoint $a$ is called \emph{quasi-regular}, if the following statement is true:
For every $\la\in\bb C$ there is a $c\in(a,b)$ such that $\Xi_\la\cap(a,c)=\emptyset$ and all balanced solutions of $Ju'+qu=\la wu$ on $(a,c)$ are in $L^2(w\chi_{(a,c)})$. 
An analogous definition is, of course, made for $b$.
We will soon see that regular endpoints are quasi-regular.
However, if $w$ is $0$ near $a$ while $q$ is not finite near $a$, then $a$ is quasi-regular but not regular. 
A singular endpoint that is not quasi-regular is called \emph{truly singular}%
\footnote{In the classical case, where $n=2$, $J=\sm{0&-1\\ 1&0}$, $q=\sm{v&0\\ 0&-1/p}$ and $w=\sm{r&0\\ 0&0}$ with $v$, $1/p$ and $r$ real and locally integrable with respect to Lebesgue measure, one uses the terms \emph{limit-circle} for quasi-regular and \emph{limit-point} for truly singular. These go back to Weyl's seminal 1910 paper \cite{41.0343.01}.}.

We can now state the following existence and uniqueness theorem for solutions of initial value problems. 
\begin{theorem}[\cite{MR5052520}, Topics 10.3.12 and 10.3.13]\label{EUIVP}
Suppose Hypothesis \ref{hyp:m} holds, $\la\in\bb C$, $f\in L^2(w)$, $x_0\in(c,d)\subset(a,b)$ and $u_0\in\bb C^n$.
If $\Xi_\la\cap(c,d)=\emptyset$ then the initial value problem $Ju'+(q-\la w)u=wf$, $u(x_0)=u_0$ has a unique balanced solution in the interval $(c,d)$.
We may also pose the initial conditions $u^+(x_0)=u_0$ or $u^-(x_0)=u_0$ to obtain the same result. 
If $c=a$ and $a$ is a regular endpoint we also get the existence and uniqueness of a solution $u$ satisfying $u^+(a)=u_0$ and an analogous statement is true if $d=b$ and $b$ is regular.
\end{theorem}

\begin{corollary}
A regular endpoint is quasi-regular.
\end{corollary}

\begin{proof}
Suppose $a$ is a regular endpoint for $(q,w)$ and $\la$ is a complex number. 
Since $q$ and $w$ are finite measures near $a$ it follows that $a$ cannot be a limit point of $\Xi_\la$, i.e., $\inf\Xi_\la>a$.
Then Theorem \ref{EUIVP} shows that any balanced solution of $Ju'+qu=\la wu$ is of bounded variation near $a$ and hence in $L^2(w\chi_{(a,b')})$ where $b'=\inf\Xi_\la$ or, if $\Xi_\la$ is empty, $b'$ is any point in $(a,b)$. 
Thus $a$ is quasi-regular.
A similar proof works, of course, for $b$. 
\end{proof}

We now consider how balanced solutions of the homogeneous equation $Ju'+{(q-\la w)u}=0$ depend on the parameter $\la\in\bb C$. 
\begin{theorem}[\cite{MR4047968}, Theorem 2.7]\label{entiresol}
Suppose Hypothesis \ref{hyp:m} holds, $f\in L^2(w)$, $x_0\in(c,d)\subset(a,b)$, $\Omega$ is a an open subset of $\bb C$, $u_0:\Omega\to\bb C^n$ is holomorphic, and $\Xi_\la\cap(c,d)=\emptyset$ for all $\la\in\Omega$.
Let $u:(c,d)\times\Omega\to\bb C^n$ be the balanced solution of $Ju'+(q-\la w)u=wf$ on $(c,d)$ satisfying $u(x_0,\la)=u_0(\la)$ for all $\la\in\Omega$. 
Then $u(x,\cdot)$ is holomorphic for every $x\in(c,d)$.
Again if $c=a$ (or $d=b$) and $a$ (or $b$) is regular the claim remains true if $x_0=a$ and $u^+(a,\la)=u_0(\la)$ (or $x_0=b$ and $u^-(b,\la)=u_0(\la)$).
\end{theorem}

We will make frequent use of Green's formula or Lagrange's identity which is the following statement:
If $Ju'+qu=wf$, $Jv'+qv=wg$, and $u$ and $v$ are balanced, it follows from the integration by parts formula \eqref{ibypt} that
\begin{equation}\label{Lagrange}
(v^*Ju)^+(d)-(v^*Ju)^-(c)=\int_{[c, d]}(v^*wf-g^*wu)
\end{equation}
whenever $[c,d]\subset(a,b)$.
This fact has the following consequences:
\begin{enumerate}
  \item If $u,f,v$, and $g$ are in $\pc L^2(w\chi_{(a,d)})$ equation \eqref{Lagrange} implies that the limit $(v^*Ju)^+(a)$ exists. 
Similarly, $(v^*Ju)^-(b)$ exists if $u,f,v$, and $g$ are in $\pc L^2(w\chi_{(c,b)})$.
\item If $f=\la u$ and $g=\ol v$, then $(v^*Ju)^-(s)=(v^*Ju)^+(t)$ for all $s,t\in(a,b)$.
\end{enumerate}

Theorems \ref{EUIVP} and \ref{entiresol} imply that we have fundamental matrices of balanced solutions for $Ju'+(q-\la w)u=0$ where each entry, as a function of $\la$, is holomorphic on open subsets of $\bb C\setminus\Lambda$.  
We are particularly interested in fundamental matrices satisfying the Wronskian relationship
\begin{equation}\label{wronski}
U^\pm(x,\ol)^* J U^\pm(x,\la)=J
\end{equation}
for all $\la\in\bb C\setminus\Lambda$ and all $x\in(a,b)$. 
The second consequence of Lagrange's identity \eqref{Lagrange} shows that equation \eqref{wronski} holds for all $x$ if it holds for one.
Specifying, for example, a symplectic initial condition $S$ for $U$ at some point $x_0\in(a,b)$, i.e., $U^+(x_0,\la)=S$ or $U^-(x_0,\la)=S$, then proves the existence of fundamental matrices satisfying equation \eqref{wronski}. 

If $U(\cdot,\la)$ satisfies \eqref{wronski} the variation of constants formula holds in the following form. 
\begin{theorem}[\cite{MR4047968}, Lemma 3.3]\label{varconst2}
For $\la\not\in\Lambda$ let $U(\cdot,\la)$ be a fundamental matrix for $Ju'+qu=\la wu$ that satisfies equation \eqref{wronski} for all $x\in(a,b)$.
Let $x_0\in(a,b)$ be a point of continuity for both $Q$ and $W$, the antiderivatives of $q$ and $w$.
If an endpoint is regular, $x_0$ may also be chosen to be that endpoint.

If $f\in \pc L^1_\loc(w)$ and if the functions $u^\pm$ are given by
\begin{equation}\label{150812.1}
u^-(x,\la)=U^-(x,\la)J^{-1}\int_{(x_0,x)} U(\cdot,\ol)^* wf, \quad x\geq x_0
\end{equation}
and
\begin{equation}\label{150812.2}
u^+(x,\la)=-U^+(x,\la)J^{-1}\int_{(x,x_0)} U(\cdot,\ol)^*wf, \quad x\leq x_0,
\end{equation}
then $u^\#$ is a balanced solution of the initial value problem $Ju'+qu=\la wu +wf$, $u(x_0)=0$.
Conversely, given a balanced solution $u$ of this initial value problem, the associated functions $u^\pm$ are given by formulas \eqref{150812.1} and \eqref{150812.2}.
\end{theorem}
The proof given in \cite{MR4047968} assumes that $U(x_0,\la)=\id$, however, it uses only that equation \eqref{wronski} holds everywhere.

\subsection{Linear relations associated with \texorpdfstring{$\boldsymbol{Ju'+qu=wf}$}{the DE}}
The differential equation $Ju'+qu=wf$ may be associated with the linear relations
$$\pc T_\mx=\{(u,f)\in \pc L^2(w)\times \pc L^2(w): u\in\bvl^\#((a,b))^n, Ju'+qu=wf\}$$
and
$$\pc T_\mn=\{(u,f)\in \pc T_\mx: \text{$\supp u$ is compact in $(a,b)$}\}.$$
These give rise to linear relations in the Hilbert space $L^2(w)\times L^2(w)$:
$$T_\mx=\{([u],[f])\in L^2(w)\times L^2(w):(u,f)\in\pc T_\mx\}$$
and
$$T_\mn=\{([u],[f])\in L^2(w)\times L^2(w):(u,f)\in\pc T_\mn\}.$$
Here, and below, we assume tacitly that $u\in\bvl^\#((a,b))^n$ when $([u],[f])\in T_\mx$.
$T_\mn$ and $T_\mx$ are called, respectively, the \emph{minimal} and \emph{maximal} relation associated with the differential equation $Ju'+qu=wf$.

The basic fact underlying the spectral theory for our equation is that $T_\mn^*=T_\mx$ and this was shown in \cite{MR4298818}, Theorem 3.1.
This result entails, in particular, that $T_\mn$ is a symmetric relation associated with $Ju'+qu=wf$.
Therefore we have von Neumann's relation
$$T_\mx=\ov{T_\mn}\oplus D_\iu\oplus D_{-\iu}$$
where $\ov{T_\mn}$ is the closure of $T_\mn$ and $D_\la=\{([u],\la[u])\in T_\mx\}$ for any $\la\in\bb C$.
The dimensions of $D_\la$ are constant as $\la$ varies in either the lower or upper half-plane.
The numbers $\dim  D_{\pm\iu}$ are called deficiency indices of $T_\mn$ and we use the notation $n_\pm=\dim D_{\pm\iu}$.
It was shown in Theorem 4.1 of \cite{MR4431055} that $n_\pm\leq n$.

Self-adjoint extensions of $T_\mn$ exist if and only if $n_+=n_-$.
In this case they are given as restrictions of $T_\mx$ by boundary conditions:
\begin{theorem}[\cite{MR4431055}, Theorem 4.4]\label{T:3.5}
Suppose Hypothesis \ref{hyp:m} holds and $n_+=n_-$. 
If $([v_1],[g_1])$, ..., $([v_{n_\pm}],[g_{n_\pm}])$ are linearly independent elements of $D_\iu\oplus D_{-\iu}$ such that
\begin{equation}\label{200820.1}
(g_k^*Jg_\ell)^-(b)-(g_k^*Jg_\ell)^+(a)=0\;\; \text{for $1\leq k,\ell\leq n_\pm$,}
\end{equation}
then
\begin{equation}\label{200821.1}
T=\{([u],[f])\in T_\mx: (g_j^*Ju)^-(b)-(g_j^*Ju)^+(a)=0\text{ for $j=1,...,n_\pm$}\}
\end{equation}
is a self-adjoint extension of $T_\mn$.

Conversely, if $T$ is a self-adjoint extension of $T_\mn$, then there are elements $([v_1],[g_1])$, ..., $([v_{n_\pm}],[g_{n_\pm}])$ of $D_\iu\oplus D_{-\iu}$ satisfying \eqref{200820.1} such that $T$ is given by \eqref{200821.1} 
\end{theorem}
We emphasize here that, on account of Lagrange's identity \eqref{Lagrange}, it is immaterial which representatives are chosen in the formulas stated in the theorem.

The equation $Ju'+qu=0$ may have nontrivial balanced solutions of norm $0$. 
Since $\norm u=0$ if and only if $wu=0$ it follows that such a solution is also a solution of $Ju'+(q-\la w)u=0$ for every $\la\in\bb C$. 
The vector space 
$$\pc L_0=\{u\in\bvl^\#((a,b))^n:Ju'+qu=0,wu=0\}$$ 
collects all such solutions (and the trivial solution). 
Our problem is called \emph{definite} if $\pc L_0=\{0\}$ and \emph{nondefinite} otherwise. 
Note that, if $([u],[f])=([v],[g])\in T_\mx$ then $Ju'+qu=wf=wg=Jv'+qv$ and hence that $u-v\in\pc L_0$.
Thus, if $\pc L_0=\{0\}$ it is harmless to identify $\pc T_\mx$ and $T_\mx$.

\section{Liouville transforms}\label{S:Liouville}
Let $\brx:(a,b)\to (\br a,\br b)$ be an increasing, bijective function and let the function $S:(a,b)\to \bb C^{n\times n}$ have the following properties:
\begin{enumerate}
   \item $S$ is of locally bounded variation,
  \item $S$ is everywhere invertible,
  \item $\det S$ is locally bounded away from $0$, and 
  \item $S^\pm(x)$ are symplectic for all $x\in(a,b)$.
\end{enumerate}
Then we define, for any function $u:(a,b)\to\bb C^n$, its \emph{Liouville transform} $\br u=\ms Lu:(\br a,\br b)\to\bb C^n$ by setting $u=S\br u\circ\brx$.
Introducing $\br S=S^{-1}$ we may write $(\ms Lu)\circ\brx=\br u\circ\brx=\br Su$. 
Two comments seem to be in order.
First, an increasing, bijective function is necessarily strictly increasing and continuous;
second, since the determinant of a symplectic matrix has modulus $1$, we have that $|\det S(x)|=1$ at all points $x$ where $S$ is continuous and thus at all but at most countably many points. 

\begin{lemma}
The composition of Liouville transforms is again a Liouville transform provided their domains line up properly. 
Moreover, Liouville transforms are invertible and the inverses so obtained are also Liouville transforms.
\end{lemma}

\begin{proof}
Let $\ms L_1$ be a Liouville transform that maps functions on $(a,b)$ to functions on $(\tilde a,\tilde b)$ and $\ms L_2$ a Liouville transform that maps functions on $(\tilde a,\tilde b)$ to functions on $(\br a,\br b)$. 
Then $u=S_1\tilde u\circ\brx_1$ and $\tilde u=S_2 \br u\circ\brx_2$. 
Thus
$$u=S_1T_2\tilde u\circ\brx$$
where $T_2=S_2\circ\brx_1$ and $\brx=\brx_2\circ\brx_1$.
Here $\brx$ is again increasing and bijective.
Moreover, $S_1T_2$ is everywhere invertible and of locally bounded variation.
We also have $\det(S_1T_2)=\det S_1 \det T_2$ is locally bounded away from $0$. 
Finally, $(S_1T_2)^\pm(x)=S_1^\pm(x) T_2^\pm(x)$ are symplectic for all $x\in(a,b)$. 
Thus $S_1T_2$ and $\brx=\brx_2\circ\brx_1$ generate the Liouville transform $\ms L_2\circ\ms L_1$. 
This proves the first claim.

For the second claim note that the inverse of an increasing, bijective function is again increasing and bijective. 
Thus $\br u=R u\circ\brx^{-1}$ where $R=S^{-1}\circ\brx^{-1}$.
Clearly, $R$ is everywhere invertible and $R^\pm=(S^\pm)^{-1}\circ\brx^{-1}$ are symplectic since the inverse of a symplectic matrix is also symplectic.
Since $\brx^{-1}$ is increasing $R$ is of locally bounded variation precisely when $S^{-1}$ is.
To show this note first that $1/\det S$ is of locally bounded variation since $\det S$ is locally bounded away from $0$. 
But $S^{-1}$ is the product of the adjugate matrix of $S$ and $1/\det S$, so it is indeed of locally bounded variation. 
Finally, since, as a function of locally bounded variation, $\det S$ is locally bounded, $\det R=(1/\det S)\circ\brx^{-1}$ is locally bounded away from $0$. 
\end{proof}
We emphasize that the given definition of a Liouville transform may be viewed as the composition of two such transforms: first $u\mapsto \tilde u =S^{-1}u=\br S u$ (with the identity function as coordinate transform) followed by $\tilde u\mapsto \br u=\tilde u\circ\brx^{-1}$ (with the identity matrix as multiplier).

Now suppose $w$ is a positive measure on $(a,b)$ with values in $\bb C^{n\times n}$ and $S:(a,b)\to\bb C^{n\times n}$ is a function of locally bounded variation.
Then $S$ is locally bounded and $S^*wS$ is again a positive, matrix-valued measure on $(a,b)$.
If $\brx:(a,b)\to (\br a,\br b)$ is an increasing, bijective function, then the push-forward measure $\br w=\brx_*(S^*wS)$ is a positive matrix-valued measure on $(\br a,\br b)$.
For $u,v\in L^2(w)$ we get now
$$\<u,v\>_w=\int_{(a,b)}(\br u\circ \brx)^* S^*wS\br v\circ \brx=\int_{(\br a,\br b)} \br u^* \brx_*(S^*wS) w \br v =\<\br u, \br v\>_{\br w}$$
where we used equation \eqref{gensubrule} to establish the second identity.
Thus we have shown the following theorem.
\begin{theorem}\label{luni}
Suppose $w$ is a positive measure defined on the Borel $\sigma$-algebra $\mc B((a,b))$ and $\ms L:u\mapsto \br u=(S^{-1}u)\circ\brx^{-1}$ is a Liouville transform.
Then $\br w=\brx_*(S^*wS)$ is a positive measure on the Borel $\sigma$-algebra $\mc B((\br a,\br b))$ and $\ms L$, viewed as an operator from $L^2(w)$ to $L^2(\br w)$, is unitary.
\end{theorem}

We are interested in the effect of Liouville transforms on first order systems of linear differential equations whose coefficients satisfy Hypothesis \ref{hyp:m}.
Thus we assume now that $q$ and $w$ are $n\times n$-matrices of distributions of order $0$ or, equivalently, local measures, such that $q$ is Hermitian and $w$ nonnegative. 
Suppose the pair $(u,f)$, where $u$ is balanced, satisfies the equation $Ju'+qu=wf$ posed on the interval $(a,b)$.
If $\ms L$ is a Liouville transform let $\br u=(\br Su)\circ\brx^{-1}$ and $\br f=(\br Sf)\circ\brx^{-1}$ (recall that we write $\br S$ for $S^{-1}$). 
We then ask under what circumstances $\br u$ is balanced and the pair $(\br u,\br f)$ satisfies an equation of the form $J\br u{}'+\br q\br u =\br w\br f$.

From the differential equation $Ju'+qu=wf$ we get 
$$B_+(x,0)u^+(x)-B_-(x,0)u^-(x)=\Delta_{w}(x)f(x).$$
Using that $u$ is balanced and that $B_++B_-=2J$ this gives
\begin{equation}\label{E:20261014.2}
u^\pm(x)=J^{-1}B_\mp(x,0)u(x)\pm\frac12 J^{-1}\Delta_{w}(x)f(x).
\end{equation}
Employing the definition $\br u^\pm\circ\brx=\br S{}^\pm u^\pm$ and requiring that $\br u$ is balanced shows then that
\begin{equation}\label{E:260914.1}
\br S(x)u(x)=\br u(\brx(x))=R(x)u(x) +\frac14 \Delta_{d\br S}(x)J^{-1}\Delta_{w}(x)f(x)
\end{equation}
where
\begin{multline*}
R(x)=\frac12\bigl(\br S{}^+(x)J^{-1}B_-(x,0)+\br S{}^-(x)J^{-1}B_+(x,0)\bigr)\\
 =\br S{}^\#(x)-\frac14\Delta_{d\br S}(x)J^{-1}\Delta_{q}(x).
\end{multline*}

Equation \eqref{E:260914.1} needs to hold for all choices of $(u,f)$. 
This is the case if and only if $\Delta_{d\br S}(x)J^{-1}\Delta_{w}(x)=0$ and $\br S(x)=R(x)$. 
As we will see later we also want $\br S{}^\#$ to be invertible.
All of these conditions are satisfied if $S$ and hence $\br S$ is continuous at $x$, since then $\br S{}^\pm(x)=\br S{}^\#(x)=\br S(x)$. 
To summarize, we call the Liouville transform $\ms L$ \emph{admissible} for $(q,w)$ if $\br S=S^{-1}$ satisfies the following three conditions for all $x\in(a,b)$ particularly at points $x$ where $\br S$ has a jump:
\begin{enumerate}
  \item $\br S(x)=\frac12 \bigl(\br S{}^+(x)J^{-1}B_-(x,0)+\br S{}^-(x)J^{-1}B_+(x,0)\bigr)$,
  \item $\Delta_{d\br S}(x)J^{-1}\Delta_{w}(x)=0$, and
  \item $\br S{}^\#(x)$ is everywhere invertible.
\end{enumerate}

We will now show that the pair $(\br u,\br f)=(\ms Lu,\ms Lf)$ satisfies the differential equation $J\br u{}'+\br q\br u =\br w\br f$ for appropriate choices of the coefficients $\br q$ and $\br w$ and that these satisfy, mutatis mutandis, Hypothesis \ref{hyp:m}.

Since our Liouville transform may be viewed as a composition of two Liouville transforms we break this program up in two steps.
First we let $\brx$ be the identity function on $(a,b)$ and assume that $S$ is such that the corresponding Liouville transform is admissible for $(q,w)$. 
Recall that $S^\#=(S^++S^-)/2$ and notice that $S^\#=S^+\br S{}^\#S^-$ is invertible. 
Let $\syc^*=J(JS^\#)^{-1}$ so that $\syc^*JS^\#=J$.
Define $\br q$ and $\br w$ by 
$$\br q=\syc^*(JS'+qS) \sas \br w=\syc^*wS.$$
The product rule for functions of locally bounded variation states that
$$u'=S^+\br u{}'+S'{\br u}{}^-=S^-\br u{}'+S'{\br u}{}^+.$$
The average of these two expressions gives that $u'=S^\# \br u '+ S'\br u$.
Therefore
$$0=\syc^*(Ju'+qu-wf)=\syc^*JS^\# \br u{}'+\syc^*(JS'+qS)\br u-\syc^*wS\br f=J\br u{}'+\br q\br u-\br w \br f.$$
It remains to show that, mutatis mutandis, $\br q$ and $\br w$ satisfy Hypothesis \ref{hyp:m}.
In view of Theorem \ref{luni} we also want that $\syc^*w S=S^*wS$.

On the set where $S^+=S^-$ we have $S=S^\#=S^\pm$ and, taking into account that $S$ is then symplectic, also $\syc=S$.
Thus, on this set, $\br q=S^*(JS'+qS)$ and $\br w=S^*wS$. 
Clearly, $S^*wS$ is nonnegative and $S^*qS$ is Hermitian. 
That the term $S^*JS'$ is also Hermitian on the set where $S^+=S^-$ follows from the fact that $S^*JS=J$ and the product rule. 

There are at most countably many points $x$ where $S^-(x)\neq S^+(x)$. 
Let $x$ be any of these points in the sequel. 
We have to show that
\begin{equation}\label{20261009.1}
\Delta_{\br q}(x)=\syc(x)^*\bigl(J\Delta_{dS}(x)+\Delta_{q}(x)S(x)\bigr)
\end{equation}
is Hermitian and that
\begin{equation}\label{20261009.2}
\Delta_{\br w}(x)=\syc(x)^*\Delta_{w}(x)S(x)=S(x)^*\Delta_{w}(x)S(x)
\end{equation}
which also shows that $\Delta_{\br w}(x)$ is nonnegative. 

\begin{proof}[Investigation of \eqref{20261009.1}]
Let $Z_\pm(x)=\br S{}^\pm(x)J^{-1}B_\mp(x,0)$.
Because $S^\pm\br S{}^\pm=\id$ an easy computation shows that
$$\Delta_q(x)=-J(S^+(x)Z_+(x)-S^-(x)Z_-(x)).$$ 
Since $\br S(x)=(Z_+(x)+Z_-(x))/2$ we find
\begin{multline*}
J\Delta_{dS}(x)+\Delta_q(x)S(x)=\bigl(J\Delta_{dS}(x)\br S(x)+\Delta_q(x)\bigr)S(x)\\
 =-JS^\#(x)\bigl(Z_+(x)-Z_-(x)\bigr)S(x).
\end{multline*}
Now recall that $\syc(x)^*JS^\#(x)=J$ and that $S(x)=\br S(x)^{-1}=2(Z_+(x)+Z_-(x))^{-1}$ to see that
$$\Delta_{\br q}(x)=-2J\bigl(Z_+(x)-Z_-(x)\bigr)\bigl(Z_+(x)+Z_-(x)\bigr)^{-1}.$$
We now use that $B_\pm(x,0)^*=-B_\mp(x,0)$ and that $\br S{}^\pm$ are symplectic to show that
$$Z_+^*JZ_+=J-\frac14 \Delta_{q}(x)J^{-1}\Delta_{q}(x)=Z_-^*JZ_-.$$
This implies 
$$(Z_++Z_-)^*J(Z_+-Z_-)+(Z_+-Z_-)^*J(Z_++Z_-)=0$$
and this gives $\Delta_{\br q}(x)=\Delta_{\br q}(x)^*$. 
\end{proof}

\begin{proof}[Investigation of \eqref{20261009.2}]
Our second claim will follow from the identity 
$$\bigl(\syc(x)^*-S(x)^*\bigr)\Delta_{w}(x)=0\;\;\text{or, equivalently,}\;\; 
 \bigl(\br S(x)^{*}\syc(x)^{*}-\id\bigr)\Delta_{w}(x)=0.$$
To show this we introduce the abbreviation $Z(x)=\frac14 \br S{}^\#(x)^{-1}\Delta_{d\br S}(x)J^{-1}\Delta_{q}(x)$ so that we have $\br S(x)^* =(\id-Z(x)^*)\br S{}^\#(x)^*$.
Since $S^\pm(x)$ are symplectic we have $JS^\pm(x)=\br S{}^\pm(x)^*J$ and hence $JS^\#(x)=\br S{}^\#(x)^*J$ which, in turn implies 
$\br S{}^\#(x)^{*}\syc(x)^{*}=\id$.
Therefore our condition becomes simply $Z(x)^*\Delta_{w}(x)=0$.
The term $Z(x)^*\Delta_w(x)$ involves the expression 
$$\Delta_{d\br S}(x)^*(\br S{}^\#(x)^*)^{-1}\Delta_w(x)=-J\br S{}^\#(x)^{-1}\Delta_{d\br S}(x)J^{-1}\Delta_w(x)$$ 
where the equality follows from Lemma \ref{L:2symps}.
Since the right-hand side is zero on account of $S$ being admissible we see that equation \eqref{20261009.2} is also satisfied.
\end{proof}

We now turn to the second step where $S=\id$ so that $u=\br u\circ\brx$. 
Let $\br q$ and $\br w$ be the push-forwards of $q$ and $w$, respectively, under the function $\brx$, i.e., 
$$\br q=\brx_*q \sas \br w=\brx_*w.$$ 
Clearly, $\br q$ is Hermitian and $\br w$ is nonnegative.
Now suppose $Ju'+qu=wf$, i.e., that $Ju'+qu-wf$ is the zero measure. 
Then
\begin{equation}\label{E:260921.1}
0=\brx_*(Ju'+qu-wf)=J\brx_*u'+(\brx_*q) u\circ\brx^{-1}-(\brx_*w) f\circ\brx^{-1}
\end{equation}
using Lemma \ref{L:A2}.
Since $u^-$ is a cumulative distribution function for the local measure $u'$, a cumulative distribution function for $\brx_*u'$ is given by 
$$y\mapsto 
\begin{cases}
(\brx_*u')([y_0,y))&\text{if $y\geq y_0$}\\ -(\brx_*u')([y,y_0))&\text{if $y<y_0$}
\end{cases}$$
where $y_0$ is a fixed point in $(\br a,\br b)$.
If $y\geq y_0$ we have
$$(\brx_*u')([y_0,y))=u'(\brx^{-1}([y_0,y))=u^-(\brx^{-1}(y))-u^-(\brx^{-1}(y_0))$$
and a similar expression if $y<y_0$.
We now set $\br u=u\circ\brx^{-1}$ and $\br f=f\circ\brx^{-1}$.
It follows that $\br u^-$ is a cumulative distribution function for the local measure $\brx_*u'$ and hence that $\br u{}'=\brx_*u'$. 
Therefore equation \eqref{E:260921.1} becomes
$$0=\brx_*(Ju'+qu-wf)=J\br u{}'+\br q \br u-\br w \br f.$$

We have proved  the following theorem.
\begin{theorem}\label{T:symde}
Suppose the pair $(q,w)$ satisfies Hypothesis \ref{hyp:m} and that $\ms L:u \mapsto \br u$ is a Liouville transform, given by $u=S\br u\circ\brx$, that is admissible for $(q,w)$. 
Let 
$$\br q=\brx_*(\syc^*(JS'+qS))\sas \br w=\brx_*(S^*wS)$$
where $\syc^*=J(JS^\#)^{-1}$.
Then $(\br q,\br w)$ also satisfies, mutatis mutandis, Hypothesis~\ref{hyp:m} and $(\br u,\br f)=(\ms Lu,\ms L f)$ satisfies $J\br u{}'+\br q \br u=\br w\br f$ whenever $(u,f)$ satisfies $Ju'+qu=wf$. 
\end{theorem}

Now define the linear relation 
$$\br{\pc T}_\mx=\{(\br u,\br f)\in \pc L^2(\br w)\times\pc L^2(\br w): \br u\in\bvl^\#(\br a,\br b)^n, J\br u{}'+\br q\br u=\br w\br f\}.$$

\begin{corollary}\label{C:symlr}
If the hypotheses of Theorem \ref{T:symde} hold, then $(u,f)\in \pc T_\mx$ if and only if $(\ms Lu,\ms Lf)\in \br{\pc T}_\mx$.
\end{corollary}

\begin{proof}
Suppose $(u,f)\in \pc T_\mx$.
Then, as we showed in the theorem, $(\br u,\br f)=(\ms Lu,\ms Lf)$ satisfies $J\br u{}'+\br q \br u=\br w\br f$.
But Theorem \ref{luni} shows that $\ms Lu$ and $\ms Lf$ are in $L^2(\br w)$, i.e., $(\ms Lu,\ms Lf)\in\br{\pc T}_\mx$.
The converse direction is proved in the same way. 
\end{proof}

In the special case where $q$ and $w$ are absolutely continuous with respect to Lebesgue measure and where $\brx$, $\brx^{-1}$ and $S$ are locally absolutely continuous functions we get from Lemma \ref{L:A3}
$$\brx'\br q\circ\brx=S^*(JS'+qS)\sas \brx'\br w\circ\brx=S^*wS.$$

\section{Boundary conditions for quasi-regular endpoints}\label{S:bc}
If $a$ is a regular endpoint a canonical choice of a fundamental matrix $U$ for $Ju'+qu=wf$ is the one that satisfies $U(a,\la)=\id$ (or some other constant symplectic matrix).
If $a$ is merely a quasi-regular endpoint such a canonical choice is generally not possible. 
In this section we will use a Liouville transform to define fundamental matrices determined by quasi-boundary values at $a$. 
Thus assume that $a$ is a quasi-regular endpoint for our coefficients $q$ and $w$ and (in this section only) that, in addition to Hypothesis \ref{hyp:m}, we have $\Delta_q(x)J^{-1}\Delta_w(x)=0$ for all $x\in(a,c)$.

Suppose $\Theta$ is a fundamental matrix for $Ju'+qu=0$ such that 
\begin{equation}\label{E:20261014.1}
\Theta^\pm(x)^*J\Theta^\pm(x)J^*=\id
\end{equation}
for all $x\in(a,c)$.
Note that the second consequence of Lagrange's identity \eqref{Lagrange} implies that equation \eqref{E:20261014.1} holds for all $x\in(a,c)$ if it holds for some $x_0\in(a,c)$. 
One might, for instance, choose $\Theta^+(x_0)=J^{-1}$.
If $J$ is unitary, i.e., if $J^*=J^{-1}$, then any $\Theta$ for which $\Theta^+$ or $\Theta^-$ assumes a symplectic initial condition at some point in $(a,c)$ will work. 

For functions $u$ from $(a,c)$ to $\bb C^n$ define the transform $\ms L: u\mapsto \cev u=S^{-1}u$ where
$$S=\Theta J^*.$$
According to equations \eqref{E:20261014.2} and \eqref{E:20261014.1} we have 
$$S(x)=B_\mp(x,0)^{-1}J \Theta^\pm(x)J^*=(\Theta^\pm(x)^*B_\mp(x,0))^{-1}$$ 
and this is true for either choice of the sign. 
Thus we introduce the abbreviation 
\begin{equation}\label{E:20261014.3}
\cev S(x)=S(x)^{-1}=\Theta^+(x)^*B_-(x,0)=\Theta^-(x)^*B_+(x,0).
\end{equation}
Note that $\cev S(x)=\Theta(x)^*J$ at any point $x$ where $\Delta_q(x)=0$.

\begin{lemma}
The map
\begin{equation}\label{cevec}
u\mapsto \cev u=\ms Lu=S^{-1}u,
\end{equation}
defined for all $\bb C^n$-valued functions on $(a,c)$, is a Liouville transform that is admissible for $(q,w)$. 
\end{lemma}

\begin{proof}
Since $\Theta J^*$ is a fundamental matrix for $Ju'+qu=0$ it is of locally of bounded variation and everywhere invertible. 
To prove that $\det S$ is bounded away from $0$ on compact subintervals of $(a,c)$ note first that, by equation \eqref{E:20261014.1} $|\det\Theta^+(x)|=1/|\det J|$ is constant. 
Since $\det B_\pm(x,0)=\det(J\pm\Delta_q(x)/2)$ and since $\Delta_q$ is bounded on the subinterval in question so is $\det B_\pm(\cdot,0)$. 
Therefore $\det S^{-1}$ is bounded above and $\det S$ is bounded below. 
That $S^+$ and $S^-$ are symplectic everywhere follows from \eqref{E:20261014.1}.
This shows that $\ms L$ is a Liouville transform. 

We now show that $\ms L$ is admissible for $(q,w)$. 
We have $\cev S{}^\pm(x)J^{-1}B_\mp(x,0)=\Theta^\pm(x)^*B_\mp(x,0)=\cev S(x)$ where the last identity is from equation \eqref{E:20261014.3}.
Thus, as required, $\cev S(x)$ is the average of these two expressions.
For property (2) notice that the differential equation evaluated at the singleton $\{x\}$ gives 
$$\Delta_{d\cev S}=\Theta^+(x)^*J-\Theta^-(x)^*J=\Theta(x)^*\Delta_q(x).$$
Multiplying on the right with $J^{-1}\Delta_w(x)$ gives the desired identity. 
Finally, since $\Theta$ is balanced and invertible the same is true for $\Theta^*$. 
Now observe that $\cev S{}^\#=\Theta^*J$.
\end{proof}

\begin{lemma}\label{L:4.2}
The Liouville transform $\ms L$ introduced above transforms the equation $Ju'+qu=wf$ posed on the interval $(a,c)$ into the equation $J\cev u{}'=\cev w\cev f$, also posed on $(a,c)$, where 
$$\cev w= S^*wS.$$
The point $a$ is a regular endpoint for $J\cev u{}'=\cev w\cev f$.
\end{lemma}

\begin{proof}
Theorem \ref{T:symde} applies here with $\brx$ being the identity map. 
Since $JS'+qS=(J\Theta'+q\Theta)J^*=0$ we have $\cev q=0$ as indicated. 
It remains to show that $\cev w$ is a finite measure near $a$. 
If $z\in\bb C^n$ let $z'=J^*z$.
Then $z^*\cev w z=(\Theta z')^*w\Theta z'$. 
Since $a$ is quasi-regular linear combinations of the columns of $\Theta$ are in $L^2(w\chi_{(a,c)})$, i.e., $z^*\cev w z$ is a finite measure on $(a,c)$ for any $z\in\bb C^n$.
\end{proof}

Since, as was just shown, $a$ is a regular point for the equation $J\cev u{}'=\cev w\cev f$, a variant of Lemma 2.3 in \cite{MR4431055} shows that the set 
$$\cev\Lambda=\{\lambda\in\bb C: \exists x\in(a,c): \det\cev B_-(x,\la)\det\cev B_+(x,\la)=0\}$$ 
intersects any ball of finite radius only in finitely many places.
Hence $\Omega=\bb C\setminus\cev\Lambda$ is an open set with at most countable many discrete punctures. 
Therefore we can find a fundamental matrix $\cev U:(a,c)\times\Omega\to\bb C^{n\times n}$ for the equation $J\cev u{}'=\la\cev w\cev u$ that satisfies the initial conditions $\cev U(a,\la)=U_0(\la)$ where $U_0(\la)$ is any invertible matrix. 
If $U_0$ is holomorphic on $\Omega$, then, according to Theorem \ref{entiresol}, $\cev U(x,\cdot)$ is also holomorphic on $\Omega$ for any $x\in(a,c)$. 
Moreover, $\cev U{}^\pm(x,\ol)^*J\cev U{}^\pm(x,\la)=J$ for all $x\in(a,c)$ if that is true for some $x\in(a,c)$.
Now define $U(x,\la)=\Theta(x)J^*\cev U(x,\la)$. 
Then  $U(\cdot,\la)$ is a fundamental matrix for $Ju'+qu=\la wu$ whenever $\la\in\Omega$, $U(x,\cdot)$ is holomorphic on $\Omega$ for any $x\in(a,c)$, and the Wronskian relationship
$$U^\pm(x,\ol)^*JU^\pm(x,\la)=J$$
is satisfied for all $\la\in\Omega$ and all $x\in(a,c)$.
 
Given a balanced solution $u$ of the differential equation $Ju'+(q-\la w)u=wf$ the limit $\cev u{}^+(a)=\lim_{x\downarrow a}\Theta(x)^*Ju(x)$ exists and is called the \emph{quasi-boundary value} for $u$ at $a$. 
We emphasize that, choosing another fundamental matrix $\Theta$ for $Ju'+qu=0$ that still satisfies equation \eqref{E:20261014.1}, gives rise to a different quasi-boundary value of $u$.

\section{The Weyl-Titchmarsh \texorpdfstring{$\bb M$}{M}-function}\label{S:WTM}
The Weyl-Titchmarsh $\bb M$-function for the self-adjoint linear relation $T$ described in Theorem \ref{T:3.5} has been constructed in \cite{MR4711083} at least for the case when $B_\pm(x,0)$ are invertible for every $x\in(a,b)$ as we will assume in this section and the next. 
In this section we will record, without proof, the most basic facts of that construction under the following hypothesis, a strengthening of Hypothesis \ref{hyp:m}.

\begin{hyp}\label{hyp:wt}
$J$ is a constant, invertible, skew-hermitian $n\times n$-matrix and $q$ and $w$ are $n\times n$-matrices whose entries are distributions of order $0$ on $(a,b)$; $q$ is Hermitian and $w$ is nonnegative.
Moreover, $\Lambda$ is a closed set of isolated points that does not intersect $\bb R$.
\end{hyp}

The function $\bb M$ is going to be defined on the complement of $\bb R\cup\Lambda$.
It is a holomorphic function for which the isolated points in $\Lambda$ are removable singularities and this allows to extend $\bb M$ to all of $\bb C\setminus\bb R$. 
This extension is a matrix-valued Nevanlinna function, i.e., a function defined on $\bb C\setminus\bb R$ with the following properties: (i) $\bb M$ is symmetric in the sense that $\bb M(\la)=\bb M(\ol)^*$, (ii) $\Im \bb M$ is non-negative in the upper half-plane, and (iii) $\bb M$ is analytic.
As a Nevanlinna function $\bb M$ has the representation
$$\bb M(\la)=A+B\la+\int_{\bb R}\Big(\frac{1}{t-\la}-\frac{t}{t^2+1}\Big)\tau(t)$$
where $A$ is hermitian, $B$ is non-negative, and $\tau$ is a non-negative matrix-valued measure called a spectral measure.
It plays a major role in spectral theory, in particular, $L^2(\tau)$ is the range of the Fourier transform $\ms F$, which is defined by
$$\ms Ff(\la)=\int_{(a,b)} U(\cdot,\ol)^*wf$$
as long as $f$ vanishes outside a compact subinterval of $(a,b)$ (and otherwise as a limit).

In \cite{MR4711083} it was assumed that the fundamental matrix $U$ used to determine $\bb M$ satisfies $U(x_0,\la)=\id$ for some $x_0\in(a,b)$ and all $\la\in\bb C\setminus\Lambda$. 
However, it is only required that the Wronskian relationship
$$U^\pm(x,\ol)^*JU^\pm(x,\la)=J$$
holds for some (and hence all) $x\in(a,b)$ and all $\la\in\bb C\setminus\Lambda$ and we assume henceforth that such a fundamental matrix is given.
$\bb M$ is intimately connected with the resolvent of $T$. 
If $R_\la$ denotes the resolvent of $T$ at $\la$, then $\bb M(\la)$ encodes the integrability and boundary conditions that $R_\la f$ has to satisfy. 
Moreover, $R_\la f$ must not have norm $0$ when $f\neq0$.
To see how this is achieved we begin with the variation of constants formula as recorded in Theorem \ref{varconst2}:
Assuming that $f$ vanishes outside a compact subinterval of $(a,b)$ and that $\la\in\bb C\setminus(\bb R\cup\Lambda)$ there is a representative $u$ of $R_\la f$ such that 
$$u^-(x)=U^-(x,\la)\Bigl(u_0+J^{-1}\int_{(a,x)} U(\cdot,\ol)^*wf\Bigr)$$
for an appropriate vector $u_0\in\bb C^n$.

To incorporate the integrability condition on the left we fix a point $c\in(a,b)$, let $N_-(\la)\subset\bb C^n$ be the set of all $\eta$ such that $u=U(\cdot,\la)\eta$ satisfies $\int_{(a,c)}u^*wu<\infty$, and let $P_-(\la)$ be the orthogonal projection onto $N_-(\la)$.
Near $a$ we have $u(x)=U(x,\la)u_0$ and this means that $u_0$ must be in $N_-(\la)$, i.e.,
\begin{equation}\label{icll}
(1-P_-(\la))u_0=0.
\end{equation}
As $u(x)=U(x,\la)\bigl(u_0+J^{-1}\ms Ff(\la)\bigr)$ near $b$ the integrability condition on the right reads
\begin{equation}\label{iclr}
(1-P_+(\la))u_0=-(1-P_+(\la))J^{-1}\ms Ff(\la).
\end{equation}
As stated in Theorem \ref{T:3.5} the boundary conditions are given by $(g^*Ju)^-(b)-(g^*Ju)^+(a)=0$, where $g=(g_1,...,g_{n_+})$ and the $g_j$ are appropriate elements of $\ran(D_\iu\oplus D_{-\iu})$.
Now define
$$A_-(\la)=-(g^*JU(\cdot,\la)P_-(\la))^+(a) \sas A_+(\la)=(g^*JU(\cdot,\la)P_+(\la))^-(b).$$
Then the boundary conditions are expressed as
\begin{equation}\label{bcl}
(A_+(\la)+A_-(\la))u_0=-A_+(\la)J^{-1}(\ms Ff)(\la).
\end{equation}
Finally, let $\id-\bb P$ be the orthogonal projection onto $N_0=\{\eta: U(\cdot,0)\eta\in\pc L_0\}$.
Then we add the requirement 
\begin{equation}\label{l0l}
(\id-\bb P)u_0=0
\end{equation}
to find the unique representative of $R_\la f$ for which $u_0\in N_0^\perp$.
Equations \eqref{icll} -- \eqref{l0l} are collected in
$$\bb F(\la) u_0=\bb H_\ell(\la)J^{-1}\ms Ff(\la)$$
where  
$$\bb F(\la)=\begin{pmatrix}
  \id-P_-(\la)\\
  \id-P_+(\la)\\
  \ms A_+(\la)+\ms A_-(\la)\\
  \id-\bb P
\end{pmatrix}
\sas
\bb H_\ell(\la)=\begin{pmatrix}
  0\\
  P_+(\la)-\id\\
  -\ms A_+(\la)\\
  0
\end{pmatrix}.$$
Since $\bb F(\la)$ has full column rank it has a left inverse $\bb F(\la)^\dag$.
Thus we compute $u_0=\bb F(\la)^\dag\bb H_\ell(\la)J^{-1}\ms Ff(\la)$.
As $u_0$ and $\ms Ff(\la)$ are in $\ran\bb P$ we may write
\begin{equation}\label{E:20261017.1}
u_0=\bb P\bb F(\la)^\dag\bb H_\ell(\la)J^{-1}\bb P\ms Ff(\la).
\end{equation}

There is also a (possibly different) representative $v$ of $R_\la f$ satisfying 
$$v^+(x)=U^+(x,\la)\Bigl(v_0-J^{-1}\int_{(x,b)} U(\cdot,\ol)^*wf\Bigr)$$
with some vector $v_0\in\bb C^n$.

Instead of \eqref{E:20261017.1} we obtain now
$$v_0=\bb P\bb F(\la)^\dag\bb H_r(\la)J^{-1}\bb P\ms Ff(\la)$$
where 
$$H_r(\la)=\begin{pmatrix}
  \id-P_-(\la)\\
  0\\
  \ms A_-(\la)\\
  0
\end{pmatrix}.$$

Since both $u$ and $v$ are representatives of $R_\la f$ their difference is an element of $\pc L_0$ and their average is also a representative of $R_\la f$. 
Now we define the Weyl-Titchmarsh function $\bb M:\bb C\setminus(\bb R\cup\Lambda)\to\bb C^{n\times n}$ by setting
$$\bb M(\la)=\bb P\bb F(\la)^\dag\bb H(\la)J^{-1}\bb P$$
where $\bb H(\la)=\frac12(\bb H_\ell(\la)+\bb H_r(\la))$.
Since $u$ is the average of $u^-$ and $u^+$ and $v$ the average of $v^-$ and $v^+$ one arrives at the following representative of $R_\la f$:
\begin{multline}\label{elf}
\frac{u(x)+v(x)}2=U(x,\la)\int \bigl(\bb M(\la)+\frac12 J^{-1}\sgn(x-\cdot)\bigr)U(\cdot,\ol)^*wf\\
     +\frac14 \bigl(U^+(x,\la)-U^-(x,\la)\bigr)J^{-1}U(x,\ol)^*\Delta_w(x)f(x).
\end{multline}

\section{Behavior of the Weyl-Titchmarsh function under Liouville transforms}\label{S:WTMLT}
We will need the following lemma:

\begin{lemma}
If $q$ and $w$ satisfy Hypothesis \ref{hyp:m} and if $\ms L:u\mapsto \br u=(S^{-1} u)\circ\brx^{-1}$ is a Liouville transform that is admissible for $(q,w)$, then
\begin{multline*}
\Lambda=\{\la\in\bb C: \exists x\in(a,b): \det B_-(x,\la)\det B_+(x,\la)=0\}\\
 =\{\la\in\bb C: \exists y\in(\br a,\br b): \det \br B_-(y,\la)\det \br B_+(y,\la)=0\}=\br\Lambda.
\end{multline*}
\end{lemma}

\begin{proof}
We assume first that $\brx$ is the identity and emulate the investigation of equation \eqref{20261009.1}. 
Define $Z_\pm(\cdot,\la)=\br S{}^\pm J^{-1}B_\mp(\cdot,\la)$ and use $S^\pm\br S{}^\pm=\id$ to get
$$\Delta_q-\la\Delta_w=B_+(\cdot,\la)-B_-(\cdot,\la)=JS^-Z_-(\cdot,\la)-JS^+Z_+(\cdot,\la).$$
We showed in Section \ref{S:Liouville} that $\br S=(Z_+(\cdot,0)+Z_-(\cdot,0))/2$ but since $\ms L$ is admissible it also true that  $\br{S}=(Z_+(\cdot,\la)+Z_-(\cdot,\la))/2$ for any $\la\in\bb C$.
Therefore 
$$J\Delta_{dS}\br S+\Delta_{q}-\la\Delta_{w}
% =\frac12J\bigl((S^+-S^-)(Z_+(\cdot,\la)+Z_-(\cdot,\la))+2 S^-Z_-(\cdot,\la)-2 S^+Z_+(\cdot,\la)\bigr)
% =\frac12J\bigl(S^+Z_+(\cdot,\la)+S^+Z_-(\cdot,\la)-S^-Z_+(\cdot,\la)-S^-Z_-(\cdot,\la) +2 S^-Z_-(\cdot,\la)-2 S^+Z_+(\cdot,\la)\bigr).
% =-\frac12J\bigl(S^+Z_+(\cdot,\la)-S^+Z_-(\cdot,\la)+S^-Z_+(\cdot,\la)-S^-Z_-(\cdot,\la)\bigr).
%  =-\frac12J(S^++S^-)(Z_+(\cdot,\la)-Z_-(\cdot,\la)).
  =-JS^\#(Z_+(\cdot,\la)-Z_-(\cdot,\la)).$$
From equations \eqref{20261009.1}, \eqref{20261009.2}, and the identity just obtained we get
\begin{multline*}
\br B_\pm(\cdot,\la)=J\pm\frac12\bigl(\Delta_{\br q}-\la\Delta_{\br w}\bigr)
 =J\pm\frac12C^*\bigl(J\Delta_{dS}\br S+\Delta_{q}-\la\Delta_{w}\bigr)S\\
 =J\mp\frac12C^*JS^\#(Z_+(\cdot,\la)-Z_-(\cdot,\la))S.
\end{multline*}
Using here $C^*JS^\#=J=J\br SS$ and $\br S=(Z_+(\cdot,\la)+Z_-(\cdot,\la))/2$ gives
$$\br B_\pm(\cdot,\la)=J\bigl(\br S \mp\frac12(Z_+(\cdot,\la)-Z_-(\cdot,\la))\bigr)S=JZ_\mp(\cdot,\la)S=(S^\mp)^*B_\pm(\cdot,\la)S.$$
Since $S$ and $S^\pm$ are invertible it follows that $\br B_+(x,\la)$ is invertible if and only if $B_+(x,\la)$ is invertible, similarly for $\br B_-(x,\la)$ and $B_-(x,\la)$. 
This proves the lemma when the coordinate transform $\brx$ is the identity.
Otherwise note, invoking Lemma \ref{L:A4}, that the coordinate transform merely moves around those points $y$ where $\br B_\pm(y,\la)$ is different from $J$.
\end{proof}

\subsection{\texorpdfstring{The case when $\boldsymbol{\br U=\ms LU}$}{Case I}}\label{SuS:6.1}
Suppose $q$ and $w$ satisfy Hypothesis \ref{hyp:wt}.
Let $\ms L:u\mapsto \br u=(S^{-1}u)\circ\brx^{-1}$ be a Liouville transform that is admissible for $(q,w)$.
Then $\ms L$ generates, as described in Section \ref{S:Liouville}, the coefficients $\br q$ and $\br w$ and these satisfy, mutatis mutandis, Hypothesis \ref{hyp:wt}.
Recall from Corollary \ref{C:symlr} that 
$$\br{\pc T}_{\mx}\circ\ms L=\ms L\circ\pc T_{\mx}$$
where we extend the domain of $\ms L$ to pairs of functions so that $\ms L(u,f)=(\ms Lu,\ms Lf)$.

Since $\ms L$ is unitary it maps functions of norm $0$ to functions of norm $0$ so that 
$$\br T_{\mx}\circ\ms L=\ms L\circ T_{\mx}.$$
By the same reasoning we see that the deficiency indices for $\br T_\mn$ are the same as those for $T_\mn$, i.e., $\br n_\pm=n_\pm$.
Now suppose the boundary conditions given by $(v_1,g_1)$, ..., $(v_{n_\pm},g_{n_\pm})$ determine a self-adjoint extension $T$ of $T_\mn$.
Define $(\br v_j,\br g_j)=\ms L(v_j,g_j)$ for $j=1,...,n_\pm$.
Then $(\br v_j,\br g_j)\in \br D_\iu\oplus\br D_{-\iu}$ and, using Lagrange's identity \eqref{Lagrange} and again that $\ms L$ is unitary,
$$(\br g_k^*J\br g_\ell)^-(b)-(\br g_k^*J\br g_\ell)^+(a)=0\;\; \text{for $1\leq k,\ell\leq n_\pm$,}$$
so that the $(\br v_j,\br g_j)$ determine a self-adjoint extension $\br T$ of $\br T_\mn$.
Indeed, in this case
$$\br T\circ\ms L=\ms L\circ T.$$

\begin{theorem}\label{T:6.2}
The Weyl-Titchmarsh functions $\bb M$ and $\br{\bb M}$ respectively associated with the relations $T$ and $\br T$ described above are identical, i.e.,
$$\br{\bb M}(\la)=\bb M(\la)$$
for all $\la\in\bb C\setminus(\bb R\cup\Lambda)$ provided $U(\cdot,\la)$ and $\br U(\cdot,\la)=\ms LU(\cdot,\la)$ are the fundamental matrices used for their construction.
\end{theorem}

\begin{proof}
Suppose $\bb M$ is constructed using the fundamental matrix $U$ for $Ju'+qu=\la wu$.
In particular, $U(x,\cdot)$ is holomorphic in $\bb C\setminus\Lambda$ for all $x\in(a,b)$ and satisfies the Wronskian relationship \eqref{wronski} for all $\la\in\bb C\setminus\Lambda$ and all $x\in(a,b)$.
Then 
$$\br U(\cdot,\la)=\ms LU(\cdot,\la)$$
is a fundamental matrix for $J\br u{}'+\br q\br u=\la \br w\br u$ such that $\br U(y,\cdot)$ is holomorphic on $\bb C\setminus\Lambda$ for all $y\in(\br a,\br b)$ and $\br U(y,\ol)^*J\br U(y,\la)=J$ for all $\la\in\bb C\setminus\Lambda$ and all $y\in(\br a,\br b)$.
Let $\br{\bb M}$ be the corresponding Weyl-Titchmarsh function for $\br T$. 

Suppose $(f,u)\in R_\la$. 
Then $(u,\la u+f)\in T$.
If $\br u=\ms Lu$ and $\br f=\ms Lf$ this gives $(\br u,\la\br u+\br f)\in\br T$. 
Thus $(\br f,\br u)\in \br R_\la$.
Now recall that a representative $u$ of $R_\la f$ is given in \eqref{elf} so that at any point $x$ where $\Delta_w(x)=0$ 
$$\br u(\brx(x))=\br U(\brx(x),\la)\int_{(a,b)}\bigl(\bb M(\la)+\frac12 J^{-1}\sgn(x-\cdot)\Bigr)\br U(\brx(\cdot),\ol)^*S^*wS\br f(\brx(\cdot)).$$
We let $y=\brx(x)$ and employ the substitution rule \eqref{gensubrule} to obtain
$$\br u(y)=\br U(y,\la)\int_{(\br a,\br b)}\bigl(\bb M(\la)+\frac12 J^{-1}\sgn(\brx^{-1}(y)-\brx^{-1}(\cdot))\bigr)\br U(\cdot,\ol)^*(\brx_*S^*wS)\br f.$$
Since $\brx^{-1}$ is strictly increasing we have $\sgn(\brx^{-1}(y)-\brx^{-1}(t))=\sgn(y-t)$.
Moreover, $\br w=\brx_*S^*wS$ and hence 
$$\br u(y)=\br U(y,\la)\int_{(\br a,\br b)}\bigl(\bb M(\la)+\frac12 J^{-1}\sgn(y-\cdot)\bigr)\br U(\cdot,\ol)^*\br w\br f.$$
At points $x$ where $\Delta_w(x)=0$ and hence $\Delta_{\br w}(y)=0$ this is a value for a representative of $\br R_\la \br f$.
Another one is given by
$$\br v(y)=\br U(y,\la)\int_{(\br a,\br b)}\bigl(\br{\bb M}(\la)+\frac12 J^{-1}\sgn(y-\cdot)\bigr)\br U(\cdot,\ol)^*\br w\br f$$
assuming again $\Delta_{\br w}(y)=0$.
The difference $\br u-\br v$ is an element of $\br{\pc L}_0$ and hence $\br{w}(\br u-\br v)=0$ and $(\bb M(\la)-\br{\bb M}(\la))\br{\ms F}f(\la)=0$.
It was shown in \cite{MR4711083} that $\ran\br{\ms F}=\ran\br{\bb P}$ and this implies that 
$\br{\bb M}(\la)-\br{\bb P}\bb M(\la)\br{\bb P}=\br{\bb P}(\br{\bb M}(\la)-\bb M(\la))\br{\bb P}=0$.
Since, as note above, $\br{\pc L}_0$ is the image of $\pc L_0$ under the Liouville transform $\ms L$ and since $\br U(\cdot,\la)=\ms LU(\cdot,\la)$ we get
$$\br N_0=\{\eta\in\bb C^n:\br U(\cdot,0)\eta\in\br{\pc L}_0\}=\{\eta\in\bb C^n:\ms LU(\cdot,0)\eta\in\br{\pc L}_0\}=N_0$$
and thus $\br{\bb P}=\bb P$.
This implies our claim.
\end{proof}

\subsection{\texorpdfstring{$\boldsymbol{\br U}$ and $\boldsymbol{U}$ are not related}{Case II}}\label{SuS:6.2}
Now suppose we have two relations $T$ and $\br T$ related by a Liouville transform $\ms L$ as described in Section \ref{SuS:6.1} and that $\bb M$ and $\br{\bb M}$ are associated Weyl-Titchmarsh functions, respectively determined by the fundamental matrices $U$ and $\br U$. 
According to Section \ref{S:WTM} these are only required to be holomorphic on $\bb C\setminus\Lambda$ and to satisfy the Wronskian relationship \eqref{wronski}. 
But $\br U$ may now be different from $\ms LU$. 

The problem can be recast in the following way: 
Let $U$ and $\tilde U$ be two fundamental matrices for the same equation $Ju'+(q-\la w)u=0$, both holomorphic on $\bb C\setminus\Lambda$ and satisfying \eqref{wronski}. 
These will satisfy
\begin{equation}\label{E:20261028.1}
\tilde U(\cdot,\la)=U(\cdot,\la)A(\la)
\end{equation}
for some matrix $A(\la)$ that is invertible for every $\la\in\bb C\setminus\Lambda$.
Each of $U$ and $\tilde U$ will determine a Weyl-Titchmarsh function, say $\bb M$ and $\tilde{\bb M}$. 
If $\tilde U=\ms L^{-1}\br U$, then according to Theorem \ref{T:6.2}, we have $\br{\bb M}=\tilde{\bb M}$ so that it suffices to study the relationship of $\bb M$ and $\tilde{\bb M}$. 

\begin{lemma}
Suppose the Weyl-Titchmarsh functions $\bb M$ and $\tilde{\bb M}$ associated with the self-adjoint linear relation $T$ are determined by the fundamental matrices $U$ and $\tilde U$, respectively, and that these are linked by equation \eqref{E:20261028.1}.
Then
$$\bb M(\la)=\bb P A(\la)\tilde{\bb M}(\la)A(\ol)^*\bb P \;\;\text{or, equivalently,}\;\; 
 \tilde{\bb M}(\la)=\tilde{\bb P} A(\la)^{-1} \bb M(\la) A(\ol)^{-1}{}^*\tilde{\bb P}.$$
\end{lemma}

\begin{proof}
First note that \eqref{wronski} implies 
$$U^\pm(x,\la)^{-1}J^{-1}(U^\pm(x,\la)^{-1})^*=J^{-1} \sas U^\pm(x,\la)J^{-1} (U^\pm(x,\la))^*=J^{-1}.$$
Applying these identities shows that $A(\la)J^{-1}A(\ol)^*=J^{-1}$. 

A representative of $R_\la f$ may be expressed in terms of $U$ and of $\tilde U$ using equation \eqref{elf} and its analogue. 
In the latter case, after replacing $\tilde U$ by $UA$, we get
\begin{multline*}
\tilde u(x)=U(x,\la)\int \bigl(A(\la)\tilde M(\la)A(\ol)^*+\frac12 J^{-1}\sgn(x-\cdot)\bigr)U(\cdot,\ol)^*wf\\
     +\frac14 \bigl(U^+(x,\la)-U^-(x,\la)\bigr)J^{-1}U(x,\ol)^*\Delta_w(x)f(x).
\end{multline*}
Subtracting from this equation \eqref{elf} gives 
$$\tilde u(x)-u(x)=U(x,\la) \bigl(A(\la)\tilde{\bb M}(\la)A(\ol)^*-\bb M(\la)\bigr)\ms Ff(\la).$$
Since this is in $\pc L_0$ and since $\ran\ms F=\ran\bb P$ we get $A(\la)\tilde{\bb M}(\la)A(\ol)^*\bb P-\bb M(\la)\in N_0$. 
Applying $\bb P$ on the left proves the lemma. 
\end{proof}

As indicated above this lemma has the following consequence for the Weyl-Titchmarsh functions of $T$ and $\br T$. 

\begin{theorem}\label{T:6.4}
Suppose $\bb M$ and $\br{\bb M}$ are, respectively, Weyl-Titchmarsh functions associated with the self-adjoint relations $T$ and $\br T$ satisfying $\br T\circ\ms L=\ms L\circ T$ for some Liouville transform $\ms L$. 
If $\bb M$ is determined by the fundamental matrix $U$ and $\br{\bb M}$ is determined by the fundamental matrix $\br U$, then 
$$\bb M(\la)=\bb PA(\la)\br{\bb M}(\la)A(\ol)^*\bb P$$
where $A(\la)$ is determined by the identity $\br U(\cdot,\la)=\ms L U(\cdot,\la)A(\la)$. 
\end{theorem}

\section{Special case: \texorpdfstring{$\boldsymbol{n=2}$, $\boldsymbol{a}$ quasi-regular}{n=2}}\label{S:n=2}
This section is devoted to the following theorem for the special case $n=2$. 

\begin{theorem}\label{T:7.1}
Suppose $J=\sm{0&-1\\ 1&0}$, $q$ and $w$ are $2\times 2$-matrices whose entries are distributions of order $0$ on $(a,b)$ with $q$ Hermitian and $w$ nonnegative, $\Delta_q(x)J^{-1}\Delta_w(x)=0$ for all $x$ near $a$, $\pc L_0=\{0\}$, $\Lambda$ is a closed set of isolated points that does not intersect $\bb R$, and the endpoint $a$ is quasi-regular.
If $T$ is a self-adjoint linear relation associated with the differential equation $Ju'+qu=wf$ and separated boundary conditions, then there is a fundamental matrix $U$ and a scalar-valued Nevanlinna function $m$ such that 
$$\bb M(\la)=\begin{pmatrix}0&1/2\\ 1/2&m(\la)\end{pmatrix}.$$
\end{theorem}

To prove this theorem we need the following lemma.

\begin{lemma}\label{L:5.2}
Let 
$$\mathbf B(\la)=\{(\ms Ff)(\la): \text{$f\in\pc L^2(w)$ and $f=0$ outside a compact subset of $(a,b)$}\}.$$
Under the hypotheses of Theorem \ref{T:7.1} we have $\mathbf B(\la)=\bb C^2$. 
\end{lemma}

\begin{proof}
Let $z\in \mathbf{B}(\la)^\perp$ so that $\int f^*w U(\cdot,\ol)z=0$.
Thus $u=U(\cdot,\ol)z$, a balanced solution of $Ju'+qu=\ol wu$, is perpendicular to a dense set in $L^2(w)$.
This shows that $wu=0$ and hence $u\in\mc L_0=\{0\}$. 
Hence $z=0$, i.e., $\mathbf{B}(\la)^\perp=\{0\}$.
\end{proof}

\begin{proof}[Proof of Theorem \ref{T:7.1}]
The separated boundary condition at the left endpoint $a$ is given by $(g_1^*Ju)^+(a)=0$ where $g_1\in\ran D_\iu\oplus D_{-\iu}$ satisfies $(g_1^*Ju)^-(b)=0$ for all $u\in\dom\pc T_\mx$, in particular, $(g_1^*Jg_1)^-(b)=0$. 
Because of equation \eqref{200820.1} we must therefore have $(g_1^*Jg_1)^+(a)=0$ as well. 
Let $\Theta$ be a fundamental matrix for $Ju'+qu=0$ satisfying $\Theta^+(x_0)=J^{-1}$ for some $x_0\in(a,b)$ and define $\cev g_1=S^{-1}g_1$ where $S=\Theta J^*$.
Recall from Section \ref{S:bc} that $S^{-1}(x)=\Theta(x)^* J$ at points where $\Delta_q(x)=0$. 
The first consequence of Lagrange's identity \eqref{Lagrange} guarantees that $\cev g_1$ has a limit $\cev g_1(a)$ at $a$. 
Employing equation \eqref{E:20261014.1} gives 
$$\cev g_1(a)^* J\cev g_1(a)=(g_1^*Jg_1)^+(a)=0.$$ 
This implies that $\cev g_1(a)$ is a multiple of a real vector and we may, in fact, choose this vector (and even $\cev g_1(a)$ itself) to be $(-\sin\alpha,\cos\alpha)^\top$ for some number $\alpha\in[0,\pi)$.

Now let $\cev w=S^*wS$ and let $\cev U(\cdot,\la)$ be a fundamental matrix for $J\cev u'=\la \cev w \cev u$ satisfying the initial condition
\begin{equation}\label{uuthph}
\cev U(a,\la)=\bsm{\cos\alpha&-\sin\alpha\\ \sin\alpha&\cos\alpha}.
\end{equation}
It was shown in Section \ref{S:bc} that $U(\cdot,\la)=S\cev U(\cdot,\la)=\Theta J^*\cev U(\cdot,\la)$ is then a fundamental matrix for $Ju'+qu=\la wu$.
We also have that $U(x,\cdot)$ is holomorphic on $\bb C\setminus\Lambda$ for all $x\in(a,b)$ and that 
$$U^\pm(x,\ol)^*JU^\pm(x,\la)=J$$
for all $\la\in\bb C\setminus\Lambda$ and all $x\in(a,b)$.
We now denote the first and second column of $U$ by $\theta$ and $\phi$, respectively. 
Since $(g_1^*J\phi(\cdot,\la))^+(a)=\cev g_1(a)^*J\cev\phi(a,\la)=0$ we see that $\phi(\cdot,\la)$ satisfies the given boundary condition at $a$. 

Let us now consider the other endpoint $b$. 
If $b$ is a quasi-regular endpoint then every solution of $Ju'+qu=\la wu$ is in $\pc L^2(w)$ and we need to pose a separated boundary condition at $b$.
If $\lambda\in\bb C\setminus(\bb R\cup\Lambda)$ this will give us a one-dimensional space of solutions of $Ju'+qu=\la wu$ that satisfy this boundary condition. 
Otherwise, if $b$ is truly singular, it follows from Lemma 3.3 in \cite{MR5139776} that there is precisely a one-dimensional space of solutions of $Ju'+qu=\la wu$ that are in $\pc L^2(w)$.
In either case this one-dimensional space cannot be spanned by $\phi(\cdot,\la)$ as this would imply that $\la$ is a nonreal eigenvalue of a self-adjoint relation which is absurd. 
Thus there is a unique number $\m(\la)\in\bb C$ such that
$$\psi(\cdot,\la)=U(\cdot,\la)\bsm{1\\ m(\la)}=\theta(\cdot,\la)+m(\la)\phi(\cdot,\la)$$
spans the space of solutions mentioned above. 
The function $\psi(\cdot,\la)$ is called the Weyl-Titchmarsh solution of $Ju'+qu=\la wu$ and $m:\bb C\setminus(\bb R\cup\Lambda)\to\bb C$ is called a (scalar) Weyl-Titchmarsh function associated with $T$ (in fact, $m$ can be extended to $\bb C\setminus\bb R$).

If $f\in\pc L^2(w)$ vanishes outside a compact subinterval of $(a,b)$ then $u=R_\la f$ has a representative given by equation \eqref{elf}. 
It must satisfy the boundary condition at $a$, i.e., $(g_1^*Ju)^+(a)=\cev g_1(a)^*J\cev u(a)=0$.
Using $\cev g_1(a)^*J=(\cos\alpha,\sin\alpha)$ and equation \eqref{uuthph} gives 
$$(1,0) \Bigl(\bb M(\la)-\frac12 J^{-1}\Bigr) (\ms Ff)(\la)=0$$
where $\bb M$ is the (matrix) Weyl-Titchmarsh function determined using the fundamental matrix $U$. 
Since, according to Lemma \ref{L:5.2}, $(\mc Ff)(\la)$ may be equal to $(1,0)^\top$ for one $f$ and $(0,1)^\top$ for another, it follows that
\begin{equation}\label{2787.2}
(1,0)\Bigl(\bb M(\la)-\frac12 J^{-1}\Bigr)=(0,0).
\end{equation}

In the interval above the support of $f$ the function $u$ must be a multiple of $\psi(\cdot,\la)$.
Denoting these multiples by $c_1$ and $c_2$ when $(\mc Ff)(\la)$ is $(1,0)^\top$ or $(0,1)^\top$, respectively, we obtain for $x$ in this interval from equation \eqref{elf} that 
\begin{equation}\label{2787.3}
\bb M(\la)+\frac12 J^{-1} = \bsm{1\\ m(\la)}(c_1,c_2).
\end{equation}
Multiplying this from the left with $(1,0)$ and using \eqref{2787.2} gives
$$(c_1,c_2)=(1,0)J^{-1}=(0,1)$$
and equation \eqref{2787.3} becomes
$$\bb M(\la)=\bsm{0&1\\ 0&m(\la)}-\frac12 J^{-1}=\bsm{0&1/2\\ 1/2&m(\la)}.$$
This concludes the proof.
\end{proof}

\appendix
\section{Local measures}\label{A:lm}
Suppose $X$ is a nonempty set and $\mc R$ a collection of subsets of $X$ with the following properties:
(i) if $A,B\in\mc R$, then so are $A\cup B$ and $A\setminus B$, (ii) if $n\mapsto A_n$ is a sequence of elements of $\mc R$ and if $A=\bigcup_{n=1}^\infty A_n$ is a subset of some $B\in\mc R$, then $A\in\mc R$, and (iii) there is a sequence $n\mapsto X_n\in\mc R$ such that $X=\bigcup_{n=1}^\infty X_n$.

Such a collection $\mc R$ is a ring of sets and has been called a \emph{$\tau$-ring} on $X$ in \cite{MR5052520}.%
\footnote{Fell and Doran \cite{MR936628} and Kelley and Srinivasan \cite{MR918770} call a collection of sets satisfying (i) and (ii) a $\delta$-ring as it will always contain countable intersections of its elements. Adding condition (iii) provides a close connection between $\mc R$ and the $\sigma$-algebra it generates, see Topic 4.2.5 in \cite{MR5052520}.}
$\tau$-rings serve as domains of \emph{local measures}, viz., a function $\mu:\mc R\to\bb C$ for which $\mu(A)=\sum_{n=1}^\infty \mu(A_n)$ whenever $n\mapsto A_n$ is a sequence of pairwise disjoint elements of $\mc R$ whose union $A$ is in $\mc R$.

We emphasize that the restriction of a local measure to an element $E$ of $\mc R$ is a complex measure on $E$. 

Now suppose $X$ is a set, $\mc R$ is $\tau$-ring on $X$, and $\mu:\mc R\to\bb C$ is a local measure.
Given another set $Y$ and a function $\brx:X\to Y$ define the set $\brx_*\mc R=\{E\subset Y: \brx^{-1}(E)\in\mc R\}$ and the function $\brx_*\mu:\brx_*\mc R\to\bb C$ that assigns to $E\in \brx_*\mc R$ the number $\mu(\brx^{-1}(E))$. 
The set $\brx_*\mc R$ is called the \emph{push-forward} of $\mc R$ and the function $\brx_*\mu$ is called the \emph{push-forward} of $\mu$. 
We then have the following result.

\begin{lemma}\label{L:A1}
If $\brx$ is bijective, then $\brx_*\mc R$ is a $\tau$-ring on $Y$ and $\brx_*\mu$ is a local measure.
Moreover, $\brx_*\mu=0$ if and only if $\mu=0$. 
\end{lemma}

\begin{proof}
Suppose $A,B\in \brx_*\mc R$, i.e., $\brx^{-1}(A),\brx^{-1}(B)\in\mc R$.
Then $\brx^{-1}(A\cup B)=\brx^{-1}(A)\cup \brx^{-1}(B)\in\mc R$ and hence $A\cup B\in \brx_*\mc R$.
One shows similarly that $A,B\in \brx_*\mc R$ implies $A\setminus B\in \brx_*\mc R$.
Next assume that $n\mapsto A_n$ is a sequence in $\brx_*\mc R$ and that $A=\bigcup_{n=1}^\infty A_n$ is included in some $B\in \brx_*\mc R$.
Then $\brx^{-1}(A)=\bigcup_{n=1}^\infty \brx^{-1}(A_n)$ where $n\mapsto \brx^{-1}(A_n)$ is a sequence in $\mc R$. 
Since $\brx^{-1}(A)\subset \brx^{-1}(B)\in\mc R$ it follows that $\brx^{-1}(A)\in\mc R$ and hence that $A\in \brx_*\mc R$. 
This establishes properties (i) and (ii) of a $\tau$-ring. 

So far we have not used that $\brx$ is bijective but will do so now.
It is always true that $E\subset \brx^{-1}(\brx(E))$.
Since $\brx$ is injective we obtain $E=\brx^{-1}(\brx(E))$.
Therefore $\brx(E)\in \brx_*\mc R$ for any $E\in\mc R$. 
In particular, if $n\mapsto X_n$ is a sequence in $\mc R$ such that $X=\bigcup_{n=1}^\infty X_n$, then $n\mapsto \brx(X_n)$ is a sequence in $\brx_*\mc R$ and $\bigcup_{n=1}^\infty \brx(X_n)=\brx(X)$.
Since $\brx$ is surjective we have $\brx(X)=Y$ and this shows that property (iii) of a $\tau$-ring holds.  
 
To show that $\brx_*\mu$ is a local measure suppose that $n\mapsto A_n$ is a sequence of pairwise disjoint elements of $\brx_*\mc R$ and that $A=\bigcup_{n=1}^\infty A_n$ is also in $\brx_*\mc R$.
Since the $\brx^{-1}(A_n)$, $n\in\bb N$, are also pairwise disjoint and since their union is $\brx^{-1}(A)$ we obtain
$$(\brx_*\mu)(A)=\mu(\brx^{-1}(A))=\mu\bigl(\bigcup_{n=1}^\infty \brx^{-1}(A_n)\bigr)
 =\sum_{n=1}^\infty \mu(\brx^{-1}(A_n))=\sum_{n=1}^\infty \mu(\brx_*(A_n))$$
which proves countable additivity.

For the final claim suppose $E=\brx^{-1}(\brx(E))\in\mc R$. 
Then $\mu(E)=(\brx_*\mu)(\brx(E))$.
Thus, if $\brx_*\mu=0$ then $\mu=0$.
The converse is also true even without using that $\brx$ is bijective. 
\end{proof}

Suppose $\mu$ (and hence $\brx_*\mu$) is a positive measure.
Then, as was shown in Topic 2.2.4 of \cite{MR5052520}, it holds that $g\in L^1(\brx_*\mu)$ if and only if $g\circ \brx\in L^1(\mu)$.
In this case one has 
\begin{equation}\label{gensubrule}
\int_Y g\, \brx_*\mu = \int_X g\circ \brx\, \mu.
\end{equation}
Associated with a local measure $\mu$ is its total variation $|\mu|$, a positive measure on the $\sigma$-algebra generated by $\mc R$. 
Note that $|\brx_*\mu|=\brx_*|\mu|$ when $\brx$ is bijective. 
Moreover, according to 5.1.15 of \cite{MR5052520}, $\mu$ has a Radon-Nikodym derivative with respect to $|\mu|$ and $\brx_*\mu$ has a Radon-Nikodym derivative with respect to $\brx_*|\mu|$. 

A function $f:X\to\bb C$ is called \emph{locally integrable} with respect to $|\mu|$ if $f\chi_A\in L^1(|\mu|)$ for all $A\in\mc R$.
In this case we define the local measure $f\mu$ on $\mc R$ by $(f\mu)(A)=\int_A f\mu$.

\begin{lemma}\label{L:A2}
Suppose $X$ and $Y$ are nonempty sets and $\mc R$ is a $\tau$-ring on $X$.
Let $\mu$ and $\nu$ be local measures on $\mc R$ and $f:X\to\bb C$ a function that is locally integrable with respect to $|\mu|$.
If $\brx:X\to Y$ is bijective, then
$$\brx_*(\mu+\nu)=\brx_*\mu+\brx_*\nu \sas \brx_*(f\mu)=(f\circ \brx^{-1}) \brx_*\mu.$$
\end{lemma}

\begin{proof}
Suppose $E\in \brx_*\mc R$. 
Then 
$$(\brx_*(\mu+\nu))(E)=(\mu+\nu)(\brx^{-1}(E))=(\brx_*\mu)(E)+(\brx_*\nu)(E).$$
This proves the first claim.

For the second assume first that $f=1$.
Let $r: X\to\bb C$ and $r_*:Y\to\bb C$ be, respectively, the Radon-Nikodym derivatives of $\mu$ with respect to $|\mu|$ and of $\brx_*\mu$ with respect to $|\brx_*\mu|=\brx_*|\mu|$.
Suppose $A\in\mc R$, then $\brx(A)\in \brx_*\mc R$ as $A=\brx^{-1}(\brx(A))$. 
Since $\chi_{A}=\chi_{\brx(A)}\circ \brx$ we obtain from equation \eqref{gensubrule} that
$$(\brx_*\mu)(\brx(A))=\int_Y \chi_{\brx(A)} r_* \brx_*|\mu| =\int_X \chi_{A}\, (r_*\circ \brx) |\mu|.$$
We also have 
$$(\brx_*\mu)(\brx(A))=\mu(A)=\int_X \chi_A  r |\mu|.$$
Since Radon-Nikodym derivatives for local measure are unique we obtain that ${r_*\circ \brx}=r$ almost everywhere with respect to $|\mu|$. 

Now assume again that $E\in \brx_*\mc R$. 
Let $g=f\circ \brx^{-1}$.
Then
$$(\brx_*(f\mu))(E)=(f\mu)(\brx^{-1}(E))=\int_X \chi_{\brx^{-1}(E)} fr |\mu|=\int_X (\chi_{E}gr_*)\circ \brx\,|\mu|.$$
Using again equation \eqref{gensubrule} gives
$$(\brx_*(f\mu))(E)=\int_Y (\chi_{E}gr_*)\, \brx_*|\mu|=\int_Y \chi_{E} g\, \brx_*\mu =(g\brx_*\mu)(E)$$
which is the second claim.
\end{proof}

In this paper we are interested in local Borel measures on an interval $(a,b)$, i.e., in local measures defined on the $\tau$-ring $\mc B'((a,b))$, the collection of all Borel sets of $(a,b)$ that are included in a compact subset of $(a,b)$. 
Such local Borel measures arise as Lebesgue-Stieltjes measures of function of locally bounded variation. 
Conversely, the cumulative distribution function of a local Borel measure is of locally bounded variation. 
Fixing $x_0\in(a,b)$ a left-continuous cumulative distribution function $M$ of the local Borel measure $\mu$ is given by
$$M(x)=
\begin{cases}\mu([x_0,x))&\text{if $x\geq x_0$,}\\
            -\mu([x,x_0))&\text{if $x<x_0$.}
\end{cases}$$
Now suppose $\brx:(a,b)\to (\br a, \br b)$ is increasing and bijective and hence continuous. 
If $y_0=\brx(x_0)$ a cumulative distribution function of $\brx_*\mu$ is given by
$$(M\circ \brx^{-1})(y)=
\begin{cases}(\brx_*\mu)([y_0,y))=\mu([x_0,\brx^{-1}(y)))&\text{if $y\geq y_0$,}\\
             -(\brx_*\mu)([y,y_0))=-\mu([\brx^{-1}(y),x_0))&\text{if $y<y_0$.}
\end{cases}$$
Since $M$ is of locally bounded variation and $\brx^{-1}$ is strictly increasing it follows that $M\circ \brx^{-1}$ is of locally bounded variation and hence that $\brx_*\mu$ is a local Borel measure.
A similar argument shows that $\brx_*\mu$ is a local Borel measure when $\brx$ is decreasing and bijective.

\begin{lemma}\label{L:A3}
Suppose $\brx:(a,b)\to (\br a, \br b)$ is monotone, bijective, and that $\brx$ and $\brx^{-1}$ are locally absolutely continuous functions.
Moreover assume that $\mu$ is absolutely continuous with respect to Lebesgue measure with Radon-Nikodym derivative $M'=(\mu/\m)$. 
Then $\brx_*\mu$ is absolutely continuous with respect to Lebesgue measure and its Radon-Nikodym derivative $(\brx_*\mu/\m)$ satisfies
$$\brx'(x)(\brx_*\mu/\m)(\brx(x))=M'(x)=(\mu/\m)(x)$$
almost everywhere with respect to Lebesgue measure.
\end{lemma}

\begin{proof}
Let $M$ be a cumulative distribution function of $\mu$.
Since $\mu=dM$ is absolutely continuous with respect to $\m$ we have that $M$ is a locally absolutely continuous function. 
As $\brx^{-1}$ is monotone we get that $M\circ\brx^{-1}$ is also locally absolutely continuous and hence that $d(M\circ\brx^{-1})=\brx_*\mu$ is absolutely continuous with respect to $\m$.
This is our first claim. 

In this situation it was shown by Serrin and Varberg \cite{MR247011} that the chain rule holds, i.e., that
\begin{equation}\label{E:20261104.1}
(M\circ\zeta)'=(M'\circ\zeta) \zeta'
\end{equation}
where $\zeta=\brx^{-1}$.
Since $\brx$ and $\zeta$ are both locally absolutely continuous we get for the same reason that $\zeta'=1/(\brx'\circ\zeta)$.
This implies that
$$\brx'(M\circ\zeta)'\circ\brx=M'$$
after composing equation \eqref{E:20261104.1} with $\brx$.
This is our second claim.
\end{proof}

\begin{lemma}\label{L:A4}
Suppose $\brx:(a,b)\to (\br a, \br b)$ is monotone and bijective.
Moreover assume that $\mu$ is a discrete local measure, i.e.,
$$\mu=\sum_{n=1}^\infty c_n \delta_{x_n}$$
where $n\mapsto x_n$ is a sequence of points in $(a,b)$, $n\mapsto c_n$ is a sequence of points in $\bb C$, and $\sum_{x_n\in K} |c_n|<\infty$ for all compact sets $K\in(a,b)$.
Then $\brx_*\mu$ is the discrete measure
$$\brx_*\mu=\sum_{n=1}^\infty c_n \delta_{\brx(x_n)}.$$
\end{lemma}

\begin{proof}
Suppose $E\in \mc B'((\br a,\br b))$. 
Then $(\brx_*\mu)(E)=\mu(\brx^{-1}(E))$ since $\delta_{x_n}(\brx^{-1}(E))=\delta_{\brx(x_n)}(E)$.
\end{proof}

\bibliographystyle{plain}
%\bibliography{c:/users/weikard/box/library/all} %office
%\bibliography{c:/users/rudi/box/library/all} %home
\end{document}